\documentclass[a4paper,11pt]{article}
\usepackage{stmaryrd}
\usepackage{amsfonts}
\usepackage{bbm}
\usepackage{float}
\usepackage{amscd}
\usepackage{mathrsfs}
\usepackage{latexsym,amssymb,amsmath,amscd,amscd,amsthm,amsxtra,xypic}
\usepackage[dvips]{graphicx}
\usepackage[utf8]{inputenc}
\usepackage[T1]{fontenc}
\usepackage{lmodern}
\usepackage{amssymb}
\usepackage[all]{xy}
\usepackage{nicefrac,mathtools,enumitem}
\usepackage{microtype}
\usepackage{lipsum}
\usepackage{mathtools}
\usepackage{xfrac}

\usepackage{geometry}
\usepackage{hyperref}

\usepackage{times}
\usepackage{amsthm}
\usepackage{amsmath}
\usepackage{amsfonts}
\usepackage{amssymb}
\usepackage{mathdots}
\usepackage{color}
\usepackage{mathrsfs}
\usepackage{stmaryrd}
\SetSymbolFont{stmry}{bold}{U}{stmry}{m}{n}

\usepackage{bm}

\usepackage{tikz-cd}
\usepackage{tikz}
\usetikzlibrary{matrix,arrows,decorations.pathmorphing}
\usepackage[all]{xy}
\usepackage{graphicx}
\usepackage{subfigure}
\usepackage[toc,page]{appendix}
\usepackage[usestackEOL]{stackengine}

\usepackage{theoremref}
\newtheorem{thm}{Theorem}[section]
\newtheorem{lem}[thm]{Lemma}

\newtheorem{pro}[thm]{Proposition}

\theoremstyle{definition}

\newtheorem{rmk}[thm]{Remark}
\newtheorem{defi}[thm]{Definition}

\newcommand{\be }{\begin{equation}}
\newcommand{\ee }{\end{equation}}
\newcommand{\id }{{\rm id}}

\newcommand{\pf}{\noindent{\bf Proof.}\ }

\newcommand {\onetwo}[1]{#1^{(1)}+#1^{(2)}}
\newcommand {\zeroone}[1]{#1^{(0)}+#1^{(1)}}

\newcommand{\h}{\mathfrak h}

\def\qed{\hfill ~\vrule height6pt width6pt depth0pt}

\newcommand{\br}[1]{   [ \cdot,    \cdot  ]   }

\newcommand{\g}{\mathfrak g}

\newcommand {\IH}{\mathbb H}
\newcommand {\IR}{\mathbb{R}}

\newcommand {\UU}{\mathfrak U}

\newcommand {\Ug}{U(\g)}
\newcommand {\fml}{[\negthinspace[\hbar]\negthinspace]}

\newcommand {\wh}{\widehat}
\newcommand {\IC}{\mathbb{C}}

\title{Stokes phenomenon and reflection equations}
\author{\small Xiaomeng Xu}
\date{}
\newcommand{\Addresses}{{
  \bigskip
  \footnotesize
\noindent \textsc{School of Mathematical Sciences \& Beijing International Center
for Mathematical Research, Peking University, Beijing 100871, China}\par\nopagebreak
  \textit{E-mail address}: \texttt{xxu@bicmr.pku.edu.cn}
}}

\begin{document}
\maketitle
\begin{abstract}
In this paper, we study the Stokes phenomenon of the cyclotomic Knizhnik-Zamolodchikov equation, and prove that its two types of Stokes matrices satisfy the Yang-Baxter and reflection equations respectively. We briefly discuss its isomonodromy deformation, and its relations with cyclotomic associators, twists and quantum symmetric pairs. 
\end{abstract}
\section{Introduction}
The Stokes phenomenon, in complex analysis, states that the asymptotic behavior of functions can differ in different angular sectors surrounding a singularity. In the papers \cite{Xu, Xu2}, such jumps of the asymptotics of solutions of Knizhnik-Zamolodchikov (KZ) equations with irregular singularities, from one sector to another, are shown to encode various structures in representation theory, including the braid relation in quantum groups, ${\rm gl}_n$-crystals, and the Gelfand-Tsetlin theory. More interestingly, the relation with representation theory brings new insights into analysis problems themselves: motivated by the Gelfand-Tsetlin theory, we derive an explicit expression of Stokes matrices of certain meromorphic linear systems of ordinary differential equations via the asymptotics of solutions of associated isomonodromy equations \cite{Xu2}.
In this paper, we deepen the relation between Stokes phenomenon and representation theory, by providing a construction of the universal solutions of reflection equations via Stokes matrices. Before stating our main theorems, we first give a brief recall of the theory of reflection equations and KZ equations.

The theory of reflection equations was initiated by Cherednik \cite{Che} in the study of factorized scattering
on the half line, and by Sklyanin in the investigation of quantum integrable models with boundary conditions \cite{Skl}. Examples of explicit universal solutions of reflection equations can be found in Kulish, Sasaki and Schwiebert \cite{KSS}. 
As the theory of quantum groups is governed by Yang-Baxter equations, the theory of quantum symmetric pair is closely related to reflection equations. In particular, Balagovic and Kolb \cite{BK}, by generalizing the construction of Bao and Wang \cite{BW0} for the quantum symmetric pair of type AIII, showed that any quantum symmetric pair in Letzter’s classification \cite{Let} gives a universal solution of reflection equation, called a universal $K$-matrix.  
A categorical framework for solutions of the reflection equation was proposed by T. tom Dieck and R. Haring-Oldenburg \cite{tD, tDH}, and the universal $K$-matrix corresponds to twisted braiding
on module categories over braided monoidal categories, see e.g., Enriquez \cite{En}, Brochier \cite{Bro0,Bro}.

The Knizhnik-Zamolodchikov (KZ) equation was discovered in the study of conformal field theory \cite{KZ}. It is a local system on the configuration space of points, whose monodromy has been studied by many people and is closely related to conformal field theory, quantum groups, representation theory of affine Lie algebras, hypergeometric functions, Hecke type algebras, geometry of cycles and so on. See e.g., Cherednik \cite{Che3}, Drinfeld \cite{Drinfeld}, Kohno \cite{Kohno}, Tsuchiya-Kanie \cite{TK}, Varchenko \cite{Varchenko}, and the reference therein. The KZ equations have various of generalizations, like allowing more general $r$-matrix form \cite{Che1} and the computation of its monodromy by Cherednik \cite{Che2}, allowing discretization, known as the quantum KZ equations by Frenkel and Reshetikhin \cite{FR}, allowing irregular singularities and compatible dynamical equations, known as generalized KZ (gKZ) equations by Felder, Markov, Tarasov and Varchenko \cite{FMTV}. Another generalization is the cyclotomic KZ equation, following Leibman \cite{Lei}, Golubeva-Leksin \cite{GL}, Enriquez-Etingof \cite[Section 4.2]{EE}, which is designed to incorporate various automorphisms on Lie algebras. When the automorphism is simply an involution, the relation between its monodromy and quantum symmetric pairs, has been studied by many authors, see e.g., Enriquez \cite{En}, De Commer-Neshveyev-Tuset-Yamashita \cite{CNTY}, to some extent, generalizing the works of Drinfeld and Kohno.

In this paper, we introduce a cyclotomic KZ equation coupled with extra irregular singularities, called a generalized cyclotomic KZ $({\rm gcKZ})$ equation, see Definition \ref{def:gcKZ}. We then study the Stokes phenomenon of the gcKZ equation around irregular singularities. In the end, we show that the Stokes matrices of the gcKZ equation gives rise to universal solutions of reflection equations. 

\subsection{Stokes matrices, Yang-Baxter and reflection equations} 
Throughout this paper, let us take the complex Lie algebra $\g={\rm gl}_n$, and take the negative transpose $\tau$ as an involution of $\g$ with spectral decomposition $\g=\frak k\oplus \frak p$,
where the fixed point Lie algebra is $\frak k={\rm so}_n$. Let $\h$ be the set of diagonal matrices, as a Cartan subalgebra of $\g$.
Let $\{e_i\}_{i\in I_\pm}$ be an orthonormal basis of $\pm1$-eigenspaces $\frak k$ and $\frak p$ with respect to the Killing form of $\g$. Set $\Omega_\frak k=\sum_{i\in I_+} e_i\otimes e_i\in \frak k\otimes \frak k$, $\Omega_\frak p=\sum_{i\in I_-} e_i\otimes e_i\in \frak p\otimes \frak p$, and denote $\Omega=\Omega_\frak k +\Omega_\frak p\in \g\otimes \g$. Furthermore, let us denote the Casimir element by $C_\frak k=\sum_{i\in I_+}e_ie_i\in U(\frak k).$

The involution $\tau$ extends to an automorphism of $U(\g)$. Let $V$ be a finite dimensional
$U(\g)\rtimes \mathbb{Z}/2\mathbb{Z}$-module, with $1\in \mathbb{Z}/2\mathbb{Z}$ encoding $\tau$, and $W$ be a finite dimensional $U(\frak k)$-module. We also denote by $\tau$ the action of
$\tau$ on $V$.
We consider respectively the two equations for a $V^{\otimes 2}$-valued function $Y(z)$ and a $W\otimes V$- valued function $F(z)$,

\begin{eqnarray}\label{eq}
&&\kappa\frac{dY}{dz}=\Big(u^{(2)}+\frac{\Omega}{z}\Big)\cdot Y,\\
&&\label{eq0}\kappa\frac{dF}{dz}=\Big(u^{(1)}+\frac{2\Omega_\frak k+C_\frak k^{(1)}}{z}\Big)\cdot F,
\end{eqnarray}
where $\kappa$ is a purely imaginary number, $u\in \h_{\rm reg}(\mathbb{R})$ the set of real regular elements in $\h$, and $u^{(1)}$ and $C_\frak k^{(1)}$ denote the action of $u$ and $C_\frak k$ on the first component of $V$ in $V\otimes V$ and $W\otimes V$ respectively.

The equation \eqref{eq} has an irregular singularity at $z =\infty$ of Poincar$\rm\acute{e}$ rank $1$. Since $\frac{u}{\kappa}$ is purely imaginary, the Stokes rays of \eqref{eq}
lie on the imaginary axis, and the Stokes sectors are the right half plane $\IH_+$ and the left half plane $\IH_-$.
Then following the theory of meromorphic linear systems (see e.g., \cite{Balser}\cite{BJL} or the Appendix A), it has a unique formal power series fundamental solution $\widehat{Y}(z)\in {\rm End}(V^{\otimes 2})$ around $z=\infty$, which under the Borel-Laplace transform gives rise to a canonical holomorphic solution $Y_+(z)$ (resp. $Y_-(z)$) in $\IH_+$ (resp. in $\IH_-$). The discontinuity of the two solutions $Y_\pm$ is measured by the {\it Stokes matrices} $S(u), S_-(u)\in{\rm End}(V^{\otimes 2})$, which are determined by
\begin{eqnarray}\label{defiSmatrix}
Y_-(z;u)=Y_+(z;u)\cdot e^{\frac{-\pi i}{\kappa} [\Omega]}S(u) , \ \ \ \ \  
Y_+(z;u)=Y_-(z;u) \cdot S_-(u) e^{\frac{\pi i}{\kappa} [\Omega]}
\end{eqnarray}
where $[\Omega]\in \h\otimes \h$ is the projection of $\Omega\in\g\otimes\g$ to $\h\otimes \h$, and the first (resp. second) identity is understood to hold in $\IH_-$
(resp. $\IH_+$) after $ Y_+$ (resp. $ Y_{-}$)
has been analytically continued counterclockwise. 

In a same way, we define the Stokes matrices
$K(u), K_-(u)\in {\rm End}(W\otimes V)$ of equation \eqref{eq0} in the two Stokes sectors $\mathbb{H}_\pm$.

\begin{thm}\label{mainthm}
For any $u\in\h_{\rm reg}(\mathbb{R})$, the two Stokes matrices $K(u)\in {\rm End}(W\otimes V)$ and $S(u)\in{\rm End}(V\otimes V)$ satisfy the Yang-Baxter and the $\tau$-twisted reflection equations
\begin{eqnarray*}\label{YBeq}
S^{12}S^{13}S^{23}&=&S^{23}S^{13}S^{12} \ \in \ {\rm End}(V^{\otimes 3}),\\
K^{01}S^{21}_\tau K^{02}S^{21}&=&S^{21}K^{02}S^{12}_\tau K^{01} \ \in \ {\rm End}(W\otimes V^{\otimes 2}).
\end{eqnarray*}
Here $S_\tau:=(\tau\otimes{\rm id})S$ and we index $W$ in $W\otimes V^{\otimes 2}$ as the $0$-th component, then $S^{ij}, K^{ij}$ denote that the two components of $S, K$ act respectively on the $i$-th and $j$-th components.
\end{thm}
As in Section \ref{section2}, the $\tau$ twist in the reflection equation naturally comes from the action of $\mathbb{Z}/2\mathbb{Z}\times \mathbb{Z}/2\mathbb{Z}$ on the ${\rm gcKZ}$ equations with $2$ variables. 
\begin{rmk}
The twisted reflection equations \cite[Equation (9.17) and Remark 9.7]{BK} are introduced to unify various reflection equations associated to different quantum symmetric pairs. The twist for the reflection equation of a given quantum symmetric pair is determined by the Dynkin data that characterises it in Araki’s classification \cite{Ara}, which (in type AI case) is related to $\tau$ by an inner automorphism of $\g$.
\end{rmk}

There are various monodromy relations among $S(u)$, $K(u)$ and the Casimir elements. See e.g., Proposition \ref{twisteq}. In particular, we have the relation $K(u)=-({\rm id}\otimes \tau)K_-(u)^{-1}\in {\rm End}(W\otimes V)$, which follows from the facts that $\tau(u)=-u$ and $({\rm id}\otimes \tau)(2\Omega_\frak k+C^{(1)}_\frak k)=2\Omega_\frak k+C^{(1)}_\frak k$.

\subsection{Variation of $u$}
Let us take a root space
decomposition $\g = \mathfrak h \oplus_{\alpha\in\Delta}\mathbb{C}e_\alpha$. For any root $\alpha\in\Delta$, set $C_\alpha=e_{\alpha}e_{-\alpha}$, and $C_{\frak k,\alpha}=\frac{1}{2}(e_{\alpha}+\tau(e_{\alpha}))e_{-\alpha}$.

\begin{thm}\label{mainthm2}
As a function of $u$, the two Stokes matrices $S(u)$ and $K(u)$ satisfy respectively
\begin{eqnarray}\label{isoeq1}
&&\kappa d_\h S(u)=\frac{1}{2}
\sum_{\alpha\in\Delta}\frac{d\alpha}{\alpha}
\left[\onetwo{C_{\alpha}},S(u)\right],\\
&&\label{isoeq}\kappa d_\h K(u)=\sum_{\alpha\in\Delta}\frac{d\alpha}{\alpha}
\left[\zeroone{C_{\frak k,\alpha}},K(u)\right],
\end{eqnarray}
where $d_\h$ is the de Rham differential on $\h$, and again we index $W$ in $W\otimes V$ as the $0$-th component.
\end{thm}
A prior the Stokes matrices are defined on the real part $u\in\h_{\rm reg}(\mathbb{R})$, but one can extend them to other points by the continuation of solutions of the above differential equation.

For a local picture, the theorem implies that the type $B$ braid group representation, as the monodromy representation $\rho(u)$ of generalized cyclotomic KZ equations (see Section \eqref{section2}), for different $u$ are equivalent. In a categorical setting, see e.g., \cite{tD, tDH, En, Bro, CNTY} for various versions, the theorem implies that the $\tau$-braided module category ${\rm Rep}(\frak k)$ over (the braided monoidal category) ${\rm Rep}(\g)$, constructed from $S(u)$ and $K_+(u)$ with different $u$ are equivalent.

For a global picture, it is interesting to study the monodromy of the equation \eqref{isoeq} with respect to $u\in\h_{\rm reg}$. Based on the quantum algebra version of the above results in Section \ref{quantumalgebra}, its monodromy should be related to the braid group actions on quantum symmetric pairs in type AI \cite{Noumi}, see \cite{MolRag, Ch}, as a cyclotomic version of the Drinfeld-Kohno theorem for Casimir equations \cite{TL}.

\subsection{Comodule algebras, cyclotomic associators and twists}\label{quantumalgebra}
In this subsection, we present a quantum algebra counterpart of the above categorical construction. Let us first recall the definition of comodule algebras, see e.g., \cite[Definition 2.7]{Kolb}.
\begin{defi}\label{comodule}
Let $(H,\Delta, R)$ be a quasi-triangular bialgebra with universal $R$-matrix $R$, an algebra involution $\tau:H\rightarrow H$ such that $(\tau\otimes \tau)(R)=R$. A right $H$-comodule algebra $B$ with coaction $\Delta_B:B\rightarrow B\otimes H$ is called quasi-triangular if there exists an invertible element $K\in B\otimes H$ such that
\begin{eqnarray*}
&(K1).& K\Delta_B(b)=({\rm id}\otimes \tau)\Delta_B(b)K,\\
&(K2).& (\Delta_B\otimes{\rm id})(K)=R^{21}_\tau K^{02} R^{12},\\
&(K3).& ({\rm id}\otimes \Delta)(K)=R^{12}K^{02}R^{12}_\tau K^{01}.
\end{eqnarray*}
Here $R_\tau:=({\rm id}\otimes \tau)(R)$, we label the tensor components of $B\otimes H\otimes H$ by $0, 1, 2$. The element $K$ is called a universal $K$-matrix for the $H$-comodule algebra $B$.
\end{defi}

Set $\UU=\Ug\fml$, $\UU_\frak k=U(\frak k)\fml$, denote by $\UU\wh{\otimes}\UU$ and $\UU_\frak k\wh{\otimes}\UU$ the
completed tensor product of $\IC\fml$--modules. 
Now let us consider the equations \eqref{eq} and \eqref{eq0}, with $\kappa={\hbar}^{-1}$ a formal parameter and $u$ rescaled by ${\hbar}^{-1}$, but for functions $Y_\hbar(z)$ and $F_\hbar(z)$ valued in $\UU_\frak k\wh{\otimes}\UU$ and $\UU \wh{\otimes}\UU$ respectively. Similar to finite dimensional case, one can study the Stokes phenomenon of these equations, and introduce canonical solutions, Stokes matrices and so on. Then Theorem \ref{mainthm} carries directly to this setting. 
\begin{thm}\label{mainthm4}
The two (quantum) Stokes matrices
\begin{eqnarray}\label{hbareq1} &&S_{\hbar}(u)\in \UU\wh{\otimes}\UU \ \ \ of \ \ \frac{dY_\hbar}{dz}=\Big(u^{(2)}+\hbar\frac{\Omega}{z}\Big)\cdot Y_\hbar, \\ \label{hbareq2} && \label{ceq1} K_{\hbar}(u)\in \UU_\frak k\wh{\otimes}\UU \ \ \ of \ \ \frac{dF_\hbar}{dz}=\Big(u^{(1)}+\hbar\frac{2\Omega_\frak k^{01}+C_\frak k^{(1)}}{z}\Big)\cdot F_\hbar,
\end{eqnarray}
satisfy the Yang-Baxter and $\tau$-twisted reflection equations. Here we index $\UU_\frak k$ as the $0$-th factor.
\end{thm}
Furthermore let us introduce the two (quantum) connection matrices $C_\hbar(u)$ and $T_\hbar(u)$, as the monodromy from $0$ to $\infty$, of the equations in \eqref{hbareq1} and $\eqref{hbareq2}$ respectively. See the Appendix for the definition of connection matrices. Then we have
\begin{thm}\label{comodalgebra}
For any $u\in\h_{\rm reg}(\mathbb{R})$, $B(u)=(\UU_\frak k, T_\hbar\Delta T_\hbar^{-1})$ is a $H(u)$-comodule algebra with the universal $K$-matrix $K_\hbar(u)$, over the bialgebra $H(u)=(\UU, C_\hbar\Delta C_\hbar^{-1}, S_{\hbar}(u))$. Here $\Delta$ is the standard
(cocommutative) coproduct.
\end{thm}
Let us give a sketch of a proof of the theorem and explain its relation with the cyclotomic associators and twists \cite{GL}\cite{En}. In particular, a precursor of Definition \ref{comodule} is the notion of a quasi-reflection algebra (QRA) over a quasi-triangular quasi-bialgebra (QTQBA) \cite[Definition 4.1]{En}.
These data satisfy in particular the same axioms with $(K1)-(K3)$ with the appearance of cyclotmoic associators. A typical example \cite[Section 4.5]{En} (see also \cite[Proposition 3.6]{CNTY}) is from the monodromy of cyclotomic ${\rm KZ}$ equations, that is $B_{KZ}=(\UU_\frak k, \Delta, e^{\pi i \hbar (\Omega_{\frak k}+C_\frak k^{(2)})}, \Psi_{KZ})$ is a QRA over the QTQBA $H_{KZ}=(\UU, \Delta, e^{\pi i \hbar \Omega}, \Phi_{KZ})$,
where $\Phi_{KZ}$ is the KZ associator and $\Psi_{KZ}$ is the monodromy from $0$ to $1$ of the reduction of the cyclotomic ${\rm KZ}$ equation with $n=2$ variables
\[\frac{dF_\hbar}{dz}=\hbar(\frac{2\Omega_{\frak k}^{01}+C_\frak k^{(1)}}{z}+\frac{\Omega^{12}}{z-1}+\frac{2\Omega_{\frak k}^{12}-\Omega^{12}}{z+1})F_\hbar,\]
where $F_\hbar$ is valued in $\UU_\frak k\wh{\otimes}\UU^{ \widehat{\otimes} 2}$. 
\begin{rmk}
To be more precise, Enriquez \cite{En} considered the semidirect product $U(\g)\rtimes \mathbb{Z}/2\mathbb{Z}\llbracket\hbar\rrbracket$ with $\mathbb{Z}/2\mathbb{Z}$ encoding the involution $\tau$. His definition of a QRA over a QTQBA avoids the $\tau$ twist, and doesn't necessarily satisfy condition $(K3)$ with the cyclotmoic associator. Here we take $U(\g)$ instead of the semidirect product, thus the conditions $(K1)$ and $(K2)$ with the associator of a QRA over a QTQBA involve the $\tau$ twist, see the convention from \cite[Section 3.2]{CNTY}. Furthermore, condition $(K3)$ is a ribbon $\tau$-braid relation \cite{Bro} as explained in \cite[Proposition 3.6]{CNTY}. By the same reason, in the following we have to take a $\tau$-twisted version of the notion of twists in \cite[Section 4.2]{En}.
\end{rmk}

A ($\tau$-twisted) twist of $H_{KZ}$ and $B_{KZ}$ is a pair $(F,G)$, where $F\in \UU\wh{\otimes}\UU$ and $G\in \UU_\frak  k\wh{\otimes}\UU$ are invertible elements. Under the twist, we get
\begin{eqnarray*}
&&H^F_{KZ}=\Big(\UU, F\Delta F^{-1}, F^{21}e^{\pi i \hbar \Omega}F^{-1}, F^{23}({\rm id}\otimes \Delta)({\rm id}\otimes \Delta)(F)\Phi_{KZ}(F^{12}(\Delta\otimes {\rm id})(F))^{-1}\Big), \\ &&B^{F,G}_{KZ}=\Big(\UU_\frak k, G\Delta G^{-1}, \tau^{(1)}(G)e^{\pi i \hbar (\Omega_{\frak k}+C_\frak k^{(2)})}G^{-1},
F^{23}({\rm id}\otimes \Delta)(G) \Psi_{KZ}(F^{12}(\Delta\otimes {\rm id})(G))^{-1}\Big),
\end{eqnarray*}
and $H^F_{KZ}$ is a QRA over the QTQBA $B^{F,G}_{KZ}$. Here $\tau^{(1)}(G)$ acts on the (first) $\UU$ component of $G$.
When the associators equal to $1$, $H^F_{KZ}$ is a quasitriangular bialgebra, and $B^{F,G}_{KZ}$ becomes a comodule algebra with the universal $K$-matrix $\tau^{(1)}(G)e^{\pi i \hbar (\Omega_{\frak k}+C_\frak k^{(2)})}G^{-1}$, that is the case of Definition \ref{comodule}. Actually, we can prove that the twist of $(H_{KZ},B_{KZ})$, under the quantum connection matrices $(F=C_\hbar(u), G=T_\hbar(u))$, is $(H(u),B(u))$. 
\begin{pro}\hfill\label{twisteq}
\begin{enumerate}
\item The quantum connection matrices $(C_\hbar(u), T_\hbar(u))$ are the twist killing the cyclotomic associators $(\Phi_{KZ}, \Psi_{KZ})$.  In particular, $B(u)$ is a reflection algebra over the quasi-triangular bialgebra $H(u)$, twist equivalent to the QRA $B_{KZ}$ over the QTQBA $H_{KZ}$.

\item We have the monodromy relations $S_{\hbar}=C_\hbar^{21}e^{\pi i \hbar \Omega}C_\hbar^{-1}$ and $K_{\hbar}=\tau^{(1)}(T_\hbar)e^{\pi i \hbar (\Omega_{\frak k}+C_\frak k^{(1)})} T_\hbar^{-1}$.
\end{enumerate}
\end{pro}
\begin{rmk}
The two sides of any identity in $(2)$ computes the monodromy along a semicircle around $0$ (in anti-cloclwise direction) and the monodromy along a semicircle around $\infty$ (in cloclwise direction) respectively, which by homotopy are same. It interprets geometrically the twists in the universal $R$ and $K$ matrices. Furthermore, the relation $K_{\hbar}=\tau^{(1)}(T_\hbar)e^{\pi i \hbar (\Omega_{\frak k}+C_\frak k^{(1)})} T_\hbar^{-1}$ implies that $\tau^{(1)}(\Delta_B(b))K_\hbar=K_\hbar \Delta_B(b)$ for any $b\in B(u)$, which is the axiom $(K1)$ of a universal $K$-matrix.  
\end{rmk}
The proof of the fact that $C_\hbar(u)$ and $T_\hbar(u)$ are the twist killing respectively the associators $\Phi_{KZ}$ and $\Psi_{KZ}$ is standard, i.e., uses the monodromy relation between certain asymptotics zones of gcKZ equations for $n=3$ and $n=2$ respectively. So we skip the proof. Here we just remark that the connection matrices between regular singularities in the ${\rm gcKZ}$ equation are the associators $\Phi_{KZ}$ and $\Psi_{KZ}$. However, if we switch to the asymptotics zones near the irregular singularities, then the irregular singularities dominate, thus the connection matrices (the KZ associators) between regular singularities in the asymptotics zones are trivial. This proposition 
is a cyclotomic analog of the construction due to Toledano Laredo \cite{TL}, where he proved that the connection matrix $C_\hbar$ (of gKZ equation for $n=2$) is a Drinfeld twist killing the KZ associator $\Phi_{KZ}$.

\begin{rmk}\label{category}
The automorphism $\tau$ on $U(\g)$ defines a braided autoequivalence of $U(\g)$-Mod in the obvious way: if $V$
is a $U(\g)$-module, then $\tau(V) = V$ as a vector space, while the module structure is given by $x\cdot_\tau v = \tau(x)v$. In particular, the result in this section defines a ribbon $\tau$-braided
right module category ${\rm Rep}(\frak k)$ over the braided monoidal category ${\rm Rep}(\g)$ in the sense of Brochier \cite{Bro}.
\end{rmk}

The theory of quantum symmetric pairs was developed by Noumi,
Sugitani and Dijkhuizen for classical Lie algebras, see e.g., \cite{Noumi, NS, NDS}, and developed by Letzter \cite{Let, Let2}
for all semisimple Lie algebras via the Drinfeld-Jimbo presentation of quantized enveloping algebras \cite{Jimbo}. As before, we focus on the case $\g={\rm gl}_n$ with the involution $\tau$, and $\frak k={\rm so}_n$ the corresponding fixed Lie subalgebra. A coideal subalgebras $B$ of $U_\hbar(\g)$ was introduced by Noumi \cite{Noumi}, see also \cite[Section 7]{Let2}, and the pair $(U_\hbar(\g), B)$ is a quantum symmetric pair of type AI (here $U_\hbar(\g)$ is the topological $\mathbb{C}\fml$-version of $U_q(\g)$, and we take the topological version of the quantum symmetric pair).
In this paper, we use Stokes phenomenon to give a transcendental construction of the pair $(H(u), B(u))$ with a universal K-matrix $K_\hbar(u)$, for any $u\in\h_{\rm reg}(\mathbb{R})$.
We expect that by a cohomological argument, for any fixed $u$ the transcendental construction $(H(u), B(u))$ and algebraic one $(U_\hbar(\g), B)$ are equivalent to each other.
\begin{rmk}
The Drinfeld-Kohno theorem \cite{Drinfeld}\cite{Kohno} states that the monodromy representation of KZ equations and the “algebraic” representations of braid group $B_n$ from universal $R$-matrices are equivalent. The equivalence between $(H(u), B(u))$ and $(U_\hbar(\g), B)$ will provide a Drinfeld-Kohno type theorem for type B braid group representation, coming from gcKZ equations and the universal $K$-matrix of algebraic quantum symmetric pairs respectively. Note that the monodromy representation of cyclotomic KZ equations and ${\rm gcKZ}$ equations are equivalent by the twist equivalence in Proposition \ref{twisteq}, it will 
be closely related to the conjecture (for the type AI case) studied by Commer, Neshveyev, Tuset and Yamashita \cite[Conjecture 4.1.]{CNTY}.
\end{rmk}

\vspace{2mm}
In this paper we use the Stokes phenomenon to construct universal $K$-matrices parameterized by $\h_{\rm reg}(\mathbb{R})$. On the one hand, as pointed by Etingof to us, it motivates an interesting question to study a dynamical analog of $K$-matrices, via the Tannakian duality for the fibre functors $F_d: {\rm Rep}(U_\hbar(g))\rightarrow A$-bimodules, where $A$ is an algebra of appropriate functions on the Cartan subalgebra $\h$, as in the dynamical $R$-matrices.
On the other hand, a theory of canonical bases and Schur duality for quantum symmetric
pairs was set up by Bao and Wang \cite{BW0, BW}. Note that the description of the WKB approximation of quantum Stokes matrices $S_\hbar(u)$ via the theory of crystals has been proposed in \cite{Xu2}. It is interesting to study the further relation between the WKB approximation of quantum Stokes matrices $S_\hbar(u)$ and $K_\hbar(u)$ (as $\hbar \rightarrow \infty$) and the canonical bases of quantum groups and quantum symmetric pairs respectively.


\section{Monodromy of generalized cyclotomic Knizhnik--Zamolodchikov equations}\label{section2}
This section computes the monodromy of ${\rm gcKZ}$ equations and proves Theorem \ref{mainthm}. In particular, Section \ref{mainsec} and \ref{Braidrep} introduce the ${\rm gcKZ}$ equation and its monodromy representation. Section \ref{solsim} introduces the canonical solutions of ${\rm gcKZ}$ equations. In the end, Section \ref{Monorep} computes the monodromy of ${\rm gcKZ}$ equations and proves Theorem \ref{mainthm}.
\subsection{Generalized cyclotomic {\rm KZ} equations}\label{mainsec}
Let $\tau$ be the involution on $\g$ fixing the Lie subalgebra $\frak k$. Let $V$ be a finite dimensional
$U(\g)\rtimes \mathbb{Z}/2\mathbb{Z}$-module, and $W$ be a finite dimensional $U(\frak k)$-module.
\begin{defi}\label{def:gcKZ}
The ${\rm gcKZ}$ equation, for a function $F(z_1,...,z_n)$ of $n$ complex variables with values in $W\otimes V^{\otimes n}$, is
\begin{eqnarray}\label{gcKZ}
\kappa\frac{\partial F}{\partial z_i}=\Big(u^{(i)}+\frac{2\Omega_\frak k^{0i}+C_\frak k^{(i)}}{z_i}+\sum_{j\ne i,j=1}^n\frac{\Omega^{ij}}{z_i-z_j}+\sum_{j\ne i,j=1}^n\frac{2\Omega_{\frak k}^{ij}-\Omega^{ij}}{z_i+ z_j}\Big)\cdot F, \ for \ i=1,...,n.
\end{eqnarray}
Here $\Omega^{ij}$ or $\Omega_{\frak k}^{ij}$ means $\Omega$ or $\Omega_\frak k$ acting on the $i$-th and $j$-th factors of $W\otimes V^{\otimes n}$ (we index $W$ as the $0$-th factor), $u^{(i)}$ and $C_{\frak k}^{(i)}$ mean $u$ and $C_\frak  k$ acting on the $i$-th $V$ factor of $W\otimes V^{\otimes n}$. 
\end{defi}

We assume henceforth that $\kappa$ is a purely imaginary number and $u\in \h_{\rm reg}(\mathbb{R})$. 
\begin{pro}
The ${\rm gcKZ}$ equation is a compatible system of differential equations over the configuration space $X_n=(\mathbb{C^\times})^n\setminus \{z_i=\pm z_j\}$. 
\end{pro}
\pf Note that when $u$ vanishes, it reduces to the cyclotomic KZ equation (see e.g., \cite[Definition 3.3]{CNTY}), which is compatible (see e.g., \cite[Lemma 3.4]{CNTY}). So it is necessary to check for different $i$ and $j$ 
\begin{eqnarray*}
&&[u^{(i)},\Omega^{ij}]-[\Omega^{ij},u^{(j)}]=0,\\
&&[u^{(i)},2\Omega_{\frak k}^{ij}-\Omega^{ij}]+[2\Omega_{\frak k}^{ij}-\Omega^{ij},u^{(j)}]=0.
\end{eqnarray*} 
Without loss of generality, we assume $i=1$ and $j=2$. The first identity follows from $[u^{(1)}+u^{(2)},\Omega]=0$. By $({\rm id}\otimes \tau)( 2\Omega_{\frak k}-\Omega)=\Omega$, the second identity is equivalent to $({\rm id}\otimes \tau)[u^{(1)}-\tau(u)^{(2)},\Omega]=0$, which follows from the relation $\tau(u)=-u$ and the first identity.
\qed

\subsection{Braid groups in type $B$ and monodromy representation of ${\rm gcKZ}$ equations}\label{Braidrep}
There is a natural action of $G_{n,2}=(\mathbb{Z}/2\mathbb{Z})^n \rtimes S_n$ on $X_n$, where the symmetric group $S_n$ acts by permutation of variables, and $i$-th $\mathbb{Z}/2\mathbb{Z}$ acts on $i$-th variable by sign permutation. The action of $\tau$ on $V$ induces an action of $(\mathbb{Z}/2\mathbb{Z})^n$ on $W\otimes V^{\otimes n}$ (which acts trivially on $W$), and $S_n$ acts on $W\otimes V^{\otimes n}$ by permuting the $V$ components. Then the ${\rm gcKZ}$
equation is $G_{n,2}$-equivariant with respect to these actions, which implies that it
also induces a system on the quotient space $X_n/G_{n,2}$. Actually, by $({\rm id}\otimes \tau)(\Omega)=2\Omega_{\frak k}-\Omega$, the last two terms of the equation's coefficients can be written in a more symmetric way as $\sum_{k=1,2}\sum_{j\ne i,j=1}^n\frac{({\rm id}\otimes \tau^k(\Omega))^{ij}}{z_i+\tau^k(z_j)}$.


Following \cite{Bri}, the fundamental group of $X_n/G_{n,2}$ is isomorphic to the braid group $B_n^1$ in type $B$ with generators $\sigma, b_1,...,b_{n-1}$ and relations
\begin{eqnarray}\label{braid}
\sigma b_i&=&b_i\sigma, \ \ \ i\ge 2,\\
\sigma b_1\sigma b_1&=&b_1\sigma b_1\sigma,\\
b_ib_j &=&b_jb_i, \ \ \ |i-j| > 1,\\
b_ib_{i+1}b_i&=&b_{i+1}b_ib_{i+1}.
\end{eqnarray}
Actually choose a base point $z =(z_1,... ,z_n)$ such that $z_i\in\IR$, $0<z_1<z_2<\cdot\cdot\cdot<z_n$, then a homomorphism is given by 
$B^1_n\rightarrow \pi_1(X_n/G_{n,2});~\sigma \mapsto the \ path \ in \ Figure \ 1, \ \ b_i\mapsto the \ path \ in \ Figure \ 2.$

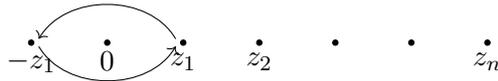
\begin{figure}[htb]
\[
  \begin{tikzpicture}
    [scale=0.5, baseline={([yshift=-.5ex]current bounding box.center)}]

    \filldraw (0,0) circle (2pt) node at (0,-0.5) {$-z_1$};
    
\filldraw (2,0) circle (2pt) node at (2,-0.5) {$0$};

    \filldraw (4,0) circle (2pt) node at (4,-0.5) {$z_1$};

    \filldraw (6,0) circle (2pt) node at (6,-0.5) {$z_{2}$};

    \filldraw (8,0) circle (2pt) node at (8.4,-0.5) {};

\filldraw (10,0) circle (2pt) node at (10.4,-0.5) {};

    \filldraw (12,0) circle (2pt) node at (12,-0.5){$z_n$};

    \path[->, bend right=60] (0.2,-0.1) edge (3.8,-0.1);
    \path[->, bend right=60] (3.8,0.1) edge (0.2,0.1);
    \end{tikzpicture}
\]
   \caption{Transposition of $z_1$ and $-z_1$ such that $z_1$ passes above $-z_1$.}
\end{figure}

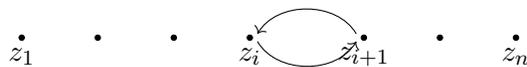
\begin{figure}[htb]
\[
  \begin{tikzpicture}
    [scale=0.5, baseline={([yshift=-.5ex]current bounding box.center)}]

    \filldraw (0,0) circle (2pt) node at (0,-0.5) {$z_1$};

    \filldraw (2,0) circle (2pt) node at (2,-0.5) {};

\filldraw (4,0) circle (2pt) node at (4,-0.5) {};

    \filldraw (6,0) circle (2pt) node at (6,-0.5) {$z_{i}$};

    \filldraw (9,0) circle (2pt) node at (9,-0.5) {$z_{i+1}$};

\filldraw (11,0) circle (2pt) node at (11,-0.5) {};

    \filldraw (13,0) circle (2pt) node at (13,-0.5){$z_n$};

    \path[->, bend right=60] (6.2,-0.1) edge (8.8,-0.1);
    \path[->, bend right=60] (8.8,0.1) edge (6.2,0.1);
    \end{tikzpicture}
\]
   \caption{Transposition of $z_i$ and $z_{i+1}$ such that $z_{i+1}$ passes above $z_{i}$.}
\end{figure}








The operators of holonomy along the paths $\sigma,b_1,...,b_{n-1}$ gives a representation of $\pi_1(X_n/G_{n,2})\cong B_n^1$. To be more precise, we fix a base point $z$, and denote by $M_0, M_i:W\otimes V^{\otimes n}\rightarrow W\otimes V^{\otimes n}$ the corresponding operator of holonomy along the paths $\sigma$ and $b_i$ in Figures $1$ and $2$. Let $s_i\in S_n$ be the permutation of the $i$-th and $i+1$-th $V$-components of $W\otimes V^{\otimes n}$. In this way, we get
\begin{pro}
The maps $\sigma\rightarrow \tau^{(1)}(M_0)$ and $b_i\mapsto  s_i(M_i)$ are a representation of the braid group $B_n^1$ in $W\otimes V^{\otimes n}$, which does
not depend on the choice of $z$ (up to isomorphism).
\end{pro}

\subsection{Canonical solutions of gcKZ equations with prescribed asymptotics}
\label{solsim} In the following, we will study the fundamental solutions of \eqref{gcKZ} with a prescribed asymptotics. To this end, let us consider the map \[P:\mathbb{C}^\times\times X_n\rightarrow X_n; \ P(z, \xi_1, . . . , \xi_n) = (z\xi_1, . . . , z\xi_n).\]
Then the pull-back of the ${\rm gcKZ}$ equation \eqref{gcKZ} under $P$ becomes 
\begin{eqnarray}\label{gcKZ1}
&&\kappa\frac{{\partial F }}{\partial z}= \Big(\sum_{i=1}^n \xi_iu^{(i)}+ 
\frac{2\sum_{0\le i<j\le n}\Omega^{ij}_{\frak k}+\sum_{i=1}^n C_{\frak k}^{(i)}}{z}\Big)\cdot F, \\ \label{gcKZ2}
&&\kappa\frac{{\partial F }}{\partial \xi_i}= \Big(z u^{(i)}+\frac{2\Omega_{\frak k}^{0i}+C_{\frak k}^{(i)}}{\xi_i}+\sum_{j\ne i,j=1}^n\frac{\Omega^{ij}}{\xi_i- \xi_j}+\sum_{j\ne i,j=1}^n\frac{2\Omega_{\frak k}^{ij}-\Omega^{ij}}{\xi_i+ \xi_j}\Big)\cdot F, \ \ \ i=1,...,n.
\end{eqnarray}
In the rest of this section, we fix a $u\in\h_{\rm reg}(\mathbb{R})$.
Let us first assume that $\xi_1,...,\xi_n$ are real, and for any $k=-1,0,1,...,n-1$, define the domain in $(\mathbb{R^\times})^n$,
\begin{eqnarray*}&&D_{-1}:=\{\xi\in \IR^n~|~\xi_1<0<\xi_2<\cdot\cdot\cdot <\xi_n\},\\
&&D_0:=\{\xi\in \IR^n~|~0<\xi_1<\xi_2<\cdot\cdot\cdot <\xi_n\},\\ &&D_k:= \{\xi\in \IR^n ~|~ 0<\xi_1< \cdot\cdot\cdot <\xi_{k-1}< \xi_{k+1}<\xi_k <\xi_{k+1}<\cdot\cdot\cdot< \xi_n\} \ for \ k=1,..., n-1.\end{eqnarray*} 

For any fixed point $\xi\in D_k$, equation \eqref{gcKZ1} becomes a meromorphic ordinary differential equation with an irregular singularity at $z=\infty$. There is a standard way to produce solutions of equation \eqref{gcKZ1} with prescribed asymptotics at $z=\infty$. 

Let us take a root space
decomposition $\g =\mathfrak h \oplus_{\alpha\in\Delta}\mathbb{C}e_\alpha$ with $(e_\alpha,e_{-\alpha})=1$, and for any root $\alpha\in\Delta$. In particular, $\h$ is the space of diagonal matrices, and for each $\alpha$ the elements $e_\alpha$, $e_{-\alpha}$ and $\tau(e_{\alpha})$ are elementary matrices taking the forms of $E_{ij}$, $E_{ji}$ and $-E_{ji}$ respectively. Set $C_\alpha=e_{\alpha}e_{-\alpha}$, $\Omega_\alpha=e_{\alpha}\otimes e_{-\alpha}$ and $C_{\frak k,\alpha}=\frac{1}{2}(e_{\alpha}+\tau(e_{\alpha}))e_{-\alpha}$, $\Omega_{\frak k,\alpha}=\frac{1}{2}(e_{\alpha}+\tau(e_{\alpha}))\otimes e_{-\alpha}$. Then

\begin{lem}\label{central}
The element
\begin{eqnarray}\label{Phi}
Y(\xi)=-\sum_{\alpha\in \Delta}\frac{1}{\alpha(u)}\Big(\sum_{1\le i<j\le n}\Big(\frac{\Omega_\alpha^{ij}}{\xi_i-\xi_j}+\frac{2\Omega_{\frak k,\alpha} ^{ij}-\Omega_\alpha^{ij}}{\xi_i+\xi_j}\Big)+\sum_{i=1}^n\frac{-2\Omega_{\frak k, \alpha}^{0i}+C_{\frak k,\alpha}^{(i)}- C_\alpha^{(i)}}{ \xi_i}\Big)
\end{eqnarray}
satisfies 
\begin{eqnarray}\label{simPhi}
[Y(\xi),\sum_{i=1}^n \xi_iu^{(i)}]=2\sum_{0\le i<j\le n}\Omega_{\frak k}^{ij}+\sum_{k=1}^n (C_\frak k^{(i)}- C_0^{(i)}),\end{eqnarray}
with $C_0:= \sum_{\alpha\in \Delta}C_\alpha^{(i)}$. Here the right hand sides of \eqref{Phi} and \eqref{simPhi} are seen as the image in the representation space ${\rm End}(W\otimes V^{\otimes n})$.
\end{lem}
\pf Note that since $\xi\in D_k$, $Y(\xi)$ is well defined, i.e., the denominator in \eqref{Phi} can not be zero. To verify \eqref{simPhi}, one uses the identities (recall that $\tau(u)=-u$) \begin{eqnarray*}&&[\sum_{k=1}^n \xi_k u^{(k)},\Omega^{ij}_\alpha]=\sum_{k=1}^n \xi_k (\delta_{ki}-\delta_{kj})\alpha(u)\Omega^{ij}_\alpha=(\xi_i-\xi_j)\alpha(u)\Omega^{ij}_\alpha,\\ 
&&[\sum_{k=1}^n \xi_k u^{(k)},\tau^{(j)}\Omega_\alpha^{ij}]=\sum_{k=1}^n  \xi_k (1-2\delta_{kj}) \tau^{(j)}[u^{(k)},\Omega_\alpha^{ij}]=(\xi_{i}+\xi_{j}) \alpha(u)\tau^{(j)}\Omega_\alpha^{ij}\\
&&[\sum_{k=1}^n \xi_k u^{(k)},-2\Omega_{\frak k, \alpha}^{0i}+C_{\frak k, \alpha}^{(i)}-C_\alpha^{(i)}]=\xi_i\alpha(u)(2\Omega_{\frak k, \alpha}^{0i}+C_{\frak k \alpha}^{i}-C_\alpha^{(i)}),
\end{eqnarray*}
and the identities 
\[({\rm id}\otimes \tau)(\Omega)=2\Omega_{\frak k}-\Omega, \hspace{5mm} \Omega_{\frak k}=\sum_{\alpha\in \Delta}\Omega_{\frak k, \alpha}, \hspace{5mm} C_{\frak k}=\sum_{\alpha\in \Delta}C_{\frak k, \alpha}. \]

\qed

\vspace{1mm}

The vector space ${\rm End}(W\otimes V^{\otimes n})$ decomposes into the direct sum of the kernel space and the image space ${\rm Im}$ of the operator ${\rm ad}_{\sum_{i=1}^n \xi_iu^{(i)}}$, and ${\rm ad}_{\sum_{i=1}^n \xi_iu^{(i)}}$ is invertible when restricts to ${\rm Im}$. Then following \ref{central}, the projection of $2\sum_{0\le i<j\le n}\Omega^{ij}_{\frak k}+\sum_{i=1}^n C_{\frak k}^{(i)}$ to the centralizer of $\sum_i \xi_iu^{(i)}$ is $\sum_{i=1}^nC_0^{(i)}$. Then we have
\begin{lem}
For any fixed $\xi\in D_k$, the ordinary differential equation \eqref{gcKZ1} has a unique formal fundamental solution taking the form \begin{eqnarray}\label{formalsum}
\widehat{F}(z;\xi)=\widehat{H}(z;\xi) e^{\frac{z}{\kappa}(\sum_i\xi_iu^{(i)})}z^{\frac{1}{\kappa}\sum_{i=1}^n C_0^{(i)}}, \ \ \ {\it for} \ \widehat{H}=1+h_1z^{-1}+h_2z^{-1}+\cdot\cdot\cdot, \end{eqnarray}
where each $h_i(\xi)\in{\rm End}(W\otimes V^{\otimes n})$.
\end{lem}
\pf
Plugging $\hat{F}$ in \eqref{formalsum} into the equation \eqref{gcKZ1} gives rise to the equation for the formal sum $\widehat{H}=1+h_1z^{-1}+h_2z^{-1}+\cdot\cdot\cdot,$
\begin{eqnarray}\label{Hhat}
\kappa\frac{d\hat{H}}{dz}+\hat{H}\cdot \Big(\sum_{i=1}^n \xi_iu^{(i)}+ 
\frac{\sum_{i=1}^n C_0^{(i)}}{z}\Big)=\Big(\sum_{i=1}^n \xi_iu^{(i)}+ 
\frac{2\sum_{0\le i<j\le n}\Omega^{ij}_{\frak k}+\sum_{i=1}^n C_{\frak k}^{(i)}}{z}\Big)\cdot \hat{H}.\end{eqnarray}
Comparing the coefficients of $z^{-1}$ on the both sides of \eqref{Hhat}, we see that $h_1(\xi)$ satisfies \begin{eqnarray}\label{simH}
[h_1(\xi),\sum_{i=1}^n \xi_iu^{(i)}]=2\sum_{0\le i<j\le n}\Omega_{\frak k}^{ij}+\sum_{k=1}^n (C_\frak k^{(i)}- C_0^{(i)}).\end{eqnarray}
Comparing the coefficients of $z^{-i}$, we see that $h_j(\xi)$, for all $j>1$, satisfies 
\begin{eqnarray}\label{simHj}
[h_j(\xi),\sum_{i=1}^n \xi_iu^{(i)}]=\Big(\kappa (j-1)+2\sum_{0\le i<j\le n}\Omega_{\frak k}^{ij}+\sum_{k=1}^n C_\frak k^{(i)}\Big)\cdot h_{j-1}(\xi)-h_{j-1}(\xi)\cdot  C_0^{(i)}.\end{eqnarray}
All $h_j(\xi)$ are then uniquely determined in a recursive way: as $j=1$, by Lemma \ref{central}, any solution $h_1(\xi)$ of \eqref{simH} can be written as a sum $Y(\xi)+X_1(\xi)$, where $X_1(\xi)$ is an element satisfying $[X_1(\xi), \sum_{i=1}^n \xi_iu^{(i)}]=0$. Different chosen $X_1(\xi)$ will make the right hand side of \eqref{simHj} for $j=2$ different. While since $\kappa$ is purely imaginary, there exists a unique $X_1(\xi)$ such that the right hand side of \eqref{simHj} for $j=2$ lives in the image space ${\rm Im}$ of the operator ${\rm ad}_{\sum_{i=1}^n \xi_iu^{(i)}}$. That is there is a unique $h_1(\xi)$ such that the relation \eqref{simHj} is well posed. Note that the operator ${\rm ad}_{\sum_{i=1}^n \xi_iu^{(i)}}$ is invertible on ${\rm Im}$, we can continue and check that each $h_j(\xi)$ is determined by \eqref{simHj} up to some term $X_j(\xi)$ in the centralizer of $\sum_{i=1}^n \xi_iu^{(i)}$, while the term $X_j(\xi)$ is further fixed by requiring the right hand side of the relation \eqref{simHj} for $j+1$ is in the space ${\rm Im}.$ It proves the existence and uniqueness of the formal fundamental solution. \qed

\vspace{2mm}
Although $\hat{H}$ is a formal power series whose radius of convergence is in general zero, it is known that (see e.g., \cite{Balser,MR}) its resummation (Borel-Laplace transform) gives a holomorphic function in each Stokes sector around $z=\infty$. In particular, by our assumption on $\kappa$ and $u$, the irregular term $\sum_{i=1}^n\frac{1}{\kappa}\xi_i u^{(i)}$ is purely imaginary, thus the only two Stokes sectors of equation \eqref{gcKZ1} are the right and left half planes $\IH_\pm$. The resummation of $\hat{H}$ in $\IH_+$ gives a unique holomorphic function $H_k: \IH_+\rightarrow {\rm End}(W\otimes V^{\otimes n})$, which can be analytically continued to the bigger sector $\wh{\IH}_+=\{\rho e^{i\theta}~|~\rho>0, -\pi<\theta<\pi\}$ and then asymptotic to $\widehat{H}$ as $z\mapsto\infty$ within $\wh{\IH}_+$ (see Appendix A for more details). Let us choose the branch of ${\rm log}(z)$, which is real on the positive real axis, with a cut along the nonnegative imaginary axis $\iota \mathbb{R}_{\ge 0}$. Thus we have 
\begin{lem}
The function $F_{D_k}:=H_k e^{\frac{z}{\kappa}(\sum_i \xi_iu^{(i)})}z^{\frac{1}{\kappa}\sum_{i=1}^n C_0^{(i)}}$ is the unique holomorphic solution of \eqref{gcKZ1} on $\IH_+$, with the prescribed asymptotics $F_{D_k}e^{-\frac{z}{\kappa}(\sum_i \xi_iu^{(i)})}z^{-\frac{1}{\kappa}\sum_{i=1}^n C_0^{(i)}}\sim \hat{H}(z)$ within $\wh{\IH}_+$. 
\end{lem} 

For any fixed $\xi$, we have found a solution of equation \eqref{gcKZ1}. Let us now consider the variation of $\xi$ in $D_k$. By the compatibility, a solution of \eqref{gcKZ1} and \eqref{gcKZ2} will take the form $F_{D_k}(z;\xi) G_k(\xi)$ for $G_k(\xi)$ being a function of $\xi$. Actually, we have

\begin{pro}\label{solution}
The function $F_{k}(z,\xi):\IH_+\times D_k\rightarrow {\rm End}(W\otimes V^{\otimes n})$, given by
\begin{eqnarray}\label{canonicalF_k}F_{k}(z,\xi):=H_k e^{\frac{1}{\kappa}(\sum_i z\xi_iu^{(i)})} z^{\sum_{i=1}^n\frac{C_0^{(i)}}{\kappa}}\cdot \prod_{i=1}^n (\xi_i)^{\frac{C_0^{(i)}}{\kappa}}\prod_{1\le i< j\le n}\Big(\frac{\xi_i-\xi_j}{\xi_i+\xi_j}\Big)^{\frac{[\Omega]^{ij}}{\kappa}}T_k,\end{eqnarray}
satisfies the equations \eqref{gcKZ1} and \eqref{gcKZ2}, where the constants $T_0=1$ and $T_k:=e^{\frac{\pi i}{\kappa}[\Omega]^{k,k+1}}$ for $k=1,...,n-1$, and $T_{-1}:=e^{\frac{\pi i}{\kappa}C_0^{(1)}}$. 
\end{pro}
\pf 
Set 
\begin{eqnarray*}&&\nabla_z=\kappa\frac{d}{dz}-\Big(\sum_{i=1}^n \xi_iu^{(i)}+ 
\frac{2\sum_{0\le i<j\le n}\Omega^{ij}_{\frak k}+\sum_{i=1}^n C_{\frak k}^{(i)}}{z}\Big),\\
&& \ \Phi \ =\sum_{i=1}^n\Big(z u^{(i)}+\frac{2\Omega_{\frak k}^{0i}+C_{\frak k}^{(i)}}{\xi_i}+\sum_{j\ne i}\frac{\Omega^{ij}}{\xi_i- \xi_j}+\sum_{j\ne i}\frac{2\Omega_{\frak k}^{ij}-\Omega^{ij}}{\xi_i+ \xi_j}\Big)d\xi_i.
\end{eqnarray*}
From the compatibility of the equations \eqref{gcKZ1} and \eqref{gcKZ2}, we have $\nabla_z(\kappa d_\xi F_k-\Phi F_k)=0$, where $d_\xi$ denotes the exterior differentiation with respect to parameters $\xi_i's$. It implies that there exists a $1$-form $B_k(\xi)$ independent of $z$ such
that $\kappa d_\xi F_k-\Phi F_k=F_k B_k$, or equivalently \begin{eqnarray}\label{B_k}
\kappa d_\xi F_k\cdot F_k^{-1}-\Phi =F_k B_k F_k^{-1}.
\end{eqnarray}
To show $B_k=0$, we compare the expansion of the both sides of this equation at $z=\infty$. 

Firstly, the formal sum $\widehat{H}=1+h_1z^{-1}+h_2z^{-2}+\cdot\cdot\cdot$ in \eqref{formalsum} satisfies \begin{eqnarray*}
\kappa\frac{dH}{dz}+H\cdot \Big(\sum_{i=1}^n \xi_iu^{(i)}+ 
\frac{\sum_{i=1}^n C_0^{(i)}}{z}\Big)=\Big(\sum_{i=1}^n \xi_iu^{(i)}+ 
\frac{2\sum_{0\le i<j\le n}\Omega^{ij}_{\frak k}+\sum_{i=1}^n C_{\frak k}^{(i)}}{z}\Big)\cdot H.\end{eqnarray*}
Comparing the coefficients of $z^{-1}$, we have seen that $h_1(\xi)$ satisfies the identity \eqref{simH}. Now differentiating \eqref{canonicalF_k} gives \[\kappa d_\xi F_k\cdot F_k^{-1}=\kappa d_\xi H_k\cdot H_k^{-1}+H_k\sum_{i=1}^n\left( u^{(i)}z+\frac{C_0^{(i)}}{\xi_i}+\frac{[\Omega]^{ij}}{\xi_i-\xi_j}-\frac{[\Omega]^{ij}}{\xi_i+\xi_j}\right)d\xi_i\cdot H_k^{-1}.\] Then by the asymptotics $H_k\sim \widehat{H}=1+h_1z^{-1}+\cdot\cdot\cdot$ in $\widehat{\IH}_+$ and the identity \eqref{simH}, we get
\begin{eqnarray}\label{sim}d_\xi F_k\cdot F_k^{-1}-\Phi\sim O(z^{-1}) \  \ at \  \ z=\infty \ \ \ in \ \ \widehat{\IH}_+.\end{eqnarray}
It thus follows from the identity \eqref{B_k} that $F_k B_k F_k^{-1}\sim O(z^{-1})$  at $z=\infty$ in $\widehat{\IH}_+$.

Secondly, by the asymptotics of $F_k$ we have 
\[F_k B_k F_k^{-1}\sim  G_k(\xi) e^{\frac{z}{\kappa}(\sum_{i} \xi_iu^{(i)})} B_k e^{-\frac{z}{\kappa}(\sum_i \xi_iu^{(i)})}G_k(\xi)^{-1} \ \ \ as \ z\rightarrow \infty \ \ in \ \widehat{\IH}_+,\]
where \[G_k(\xi)=\prod_{i=1}^n \xi_i^{\frac{C_0^{(i)}}{\kappa}}\prod_{1\le i< j\le n}\Big(\frac{\xi_i-\xi_j}{\xi_i+\xi_j}\Big)^{\frac{1}{\kappa}[\Omega]^{ij}}\cdot T_k.\]
Since the exponentials dominate, $B_k$ must commute with $e^{-\frac{z}{\kappa}(\sum_i \xi_iu^{(i)})}$, otherwise some entries of \[e^{\frac{z}{\kappa}(\sum_{i} \xi_iu^{(i)})} B_k e^{-\frac{z}{\kappa}(\sum_i \xi_iu^{(i)})}\in {\rm End}(W\otimes V^{\otimes n})\] would grow exponentially for the opening of the sector $\widehat{\IH}_+$ is
larger than $\pi$, which would contradict with the fact $F_kB_kF_k^{-1}\sim O(z^{-1})$ in $\widehat{\IH}_+$. Thus 
\[F_k B_k  F_k^{-1}\sim G_k(\xi)B_kG_k(\xi)^{-1} +O(z^{-1})  \  \ at \  \ z=\infty \ \ \ in \ \ \widehat{\IH}_+.\]
Comparing it with $F_kB_kF_k^{-1}\sim O(z^{-1})$, we obtain that $B_k=0$. It proves that $d_\xi F_k(z,\xi)=\Phi F_k(z,\xi)$. \qed

\vspace{1mm}
Proposition \ref{solution} enables us to construct unique (therefore canonical) solutions with the prescribed asymptotics at infinity, associated to the domains 
\begin{eqnarray*}&&D_{-1}=\{z\in \IR^n~|~z_1<0<z_2<\cdot\cdot\cdot <z_n\},\\
&&D_0=\{z\in \IR^n~|~0<z_1<z_2<\cdot\cdot\cdot <z_n\},\\ &&D_k=\{z~|~0<z_1<\cdot\cdot\cdot<z_{k-1}<z_{k+1}<z_k<\cdot\cdot\cdot z_n\} \ for \ k=1,..., n-1.\end{eqnarray*} 

\begin{defi}
For any $k=-1,0,...,n-1$, we call $F_k$ the canonical solution of \eqref{gcKZ}.
\end{defi}

\subsection{Monodromy representations}\label{Monorep}
To compute the monodromy of equation \eqref{gcKZ}, we take an infinite base point $z=(z_1,...,z_n)$ in $D_0$, i.e., $z_i\ll z_{i+1}$ for any $i$ (See e.g., \cite{EFK} Section 8.4 for the monodromy representation with respect to an infinite base point). Then the induced braid group representation is
\[\pi_1(X_n/G_{n,2})\rightarrow {\rm End}(W\otimes V^{\otimes n}) \ ; \ \sigma\mapsto \tau^{(1)}(F_{-1}(z)F_0(z)^{-1}), \  b_i\mapsto s_i(F_i(z) \cdot F_0(z)^{-1}),\]
where $\sigma$ and $b_i's$ are the generators of $B_n^1$, and the image is the ratio of the canonical solutions $F_0$ and $F_i$ are taken in $D_i$ (after $F_0$ has been analytic continued to $D_i$ along the path $\tau$ or $b_i$).
\begin{thm}\label{monodromrep}
The monodromy representation of the ${\rm gcKZ}_n$ equation, with respect to the infinite base point, is given by
$$\pi_1(X_n/G_{n,2})\cong B_n\rightarrow {\rm End}(W\otimes V^{\otimes n}); \ \sigma\mapsto \tau^{(1)}(K(u)), \ b_i\mapsto s_i(S(u)^{i,i+1}),$$
where $S(u)$ and $K(u)$ are the Stokes matrices of equation \eqref{eq} and \eqref{eq0} (associated to $V$ and $W$) respectively.
\end{thm}
\pf
To get the monodromy along the path $b_i$,
we need to compute the ratio of the canonical solutions $F_0$ and $F_i$. To this end, let us first find the approximation of $F_0$ in the asymptotic zone \[z_{j+1}-z_j\gg 0 \ \text{for} \ \forall j, \ \ \ \text{and} \ \ \ |\frac{z_i-z_{i+1}}{z_i-z_{j}}|, \ \ |\frac{z_i-z_{i+1}}{z_{i+1}-z_{j}}|\ll 1, \ \text{for} \ \forall j\ne i, i+1.\] 

Let us consider the equation with two variables $z_i$ and $z_{i+1}$ for a ${\rm End}(W\otimes V^{\otimes n})$ valued function
\begin{eqnarray*}
&&\kappa\frac{\partial F}{\partial z_i}= (u^{(i)}+ \frac{\Omega_{i,i+1}}{z_i-z_{i+1}})\cdot F,\\
&&\kappa\frac{\partial F}{\partial z_{i+1}}= (u^{(i+1)}+ \frac{\Omega_{i,i+1}}{z_{i+1}-z_i})\cdot F.
\end{eqnarray*}
If we let $Y_\pm(z)\in {\rm End}(W\otimes V^{\otimes n})$ be the canonical solutions of $\kappa\frac{dY}{dz}=(u^{(i+1)}+\frac{\Omega_{i,i+1}}{z})\cdot Y$ in the two Stokes sectors $\IH_\pm$, then \begin{eqnarray}\label{reduce1}
&&F^\circ_+(z_i,z_{i+1}):=Y_+(z_{i+1}-z_i)e^{\frac{z_{i}}{\kappa}(u^{(i)}+u^{({i+1})})},\\
\label{reduce2}&&F^\circ_-(z_i,z_{i+1}):=Y_-(z_{i+1}-z_{i})e^{\frac{z_{i}}{\kappa}(u^{(i)}+u^{({i+1})})}e^{\frac{\pi i}{\kappa}[\Omega]^{i,i+1}},
\end{eqnarray}
satisfy the above two equations, and \begin{eqnarray}\label{H+0}
F^\circ_\pm(z_i,z_{i+1})=H^\circ_\pm(z_i,z_{i+1}) e^{\frac{1}{\kappa}(z_iu^{(i)}+z_{i+1}u^{(i+1)})} (z_i-z_{i+1})^{\frac{1}{\kappa}[\Omega]^{i,i+1}},
\end{eqnarray}
for a function $H^\circ_+$ (resp. $H^\circ_-$) tending to $1$ as $z_{i+1}-z_i\rightarrow \infty$ (resp. $z_{i}-z_{i+1}\rightarrow \infty$). 

Now we want to compare $F^\circ_+$ and $F_0$ in the asymptotic zone.
For that let us denote 
\begin{eqnarray}\label{Ni}
N_{i,i+1}:= e^{\frac{1}{\kappa}(\sum_{j\ne i,i+1} z\xi_ju^{(j)})} \prod_{j=1}^n (z_j)^{\frac{C_0^{(j)}}{\kappa}},
\end{eqnarray}
and consider $J=F_0\cdot (F^\circ_+N_{i,i+1})^{-1}$. Using the identity \eqref{canonicalF_k} for $k=0$, the identity \eqref{H+0} and the definition \eqref{Ni} of $N_{i,i+1}$, the function $J$ can be rewritten as $J=H_0 Z_{i,i+1} {H_+^\circ}^{-1}$ in terms of $H_{0}$ and $H_+^\circ$, where \begin{eqnarray*}  
Z_{i,i+1}(z):= (z_i-z_{i+1})^{-\frac{1}{\kappa}[\Omega]^{i,i+1}}\prod_{1\le j<k\le n}\Big(\frac{z_j-z_k}{z_j+z_k}\Big)^{\frac{1}{\kappa}[\Omega]^{jk}}.
\end{eqnarray*}
Let $W_{i,i+1}:=Z_{i,i+1}(z)|_{z_{i+1}=z_i}$ be the function obtained by letting the variable $z_{i+1}$ to be $z_i$.
Since $F_+^\circ$ (thus $H_+^\circ$) commutes with $[\Omega]^{(ji)}+[\Omega]^{(j,i+1)}$ for any $j\ne i, i+1$, the function $W_{i,i+1}$ commutes with $H_+^\circ$. Then we have

\begin{eqnarray}J=H_0 Z_{i,i+1} {H_+^\circ}^{-1}=\Big(H_0 \frac{Z_{i,i+1}}{W_{i,i+1}} {H_0^\circ}^{-1} \frac{W_{i,i+1}}{Z_{i,i+1}}\Big)\cdot  Z_{i,i+1}
\end{eqnarray}
Note that all the three functions $H_0$, $H_0^\circ$ and $ \frac{Z_{i,i+1}}{W_{i,i+1}}$ asymptotically equal to $1$ in the above chosen asymptotics zone, we obtain that $J\sim Z_{i,i+1}$. Therefore, $F_0(z)\sim Z_{i,i+1} F_+^\circ(z_i,z_{i+1}) N_{i,i+1}$. Similarly, we have $F_i(z)\sim Z_{i,i+1} F_-^\circ(z_i,z_{i+1}) N_{i,i+1}$ in a different asymptotic zone obtained by replacing $z_{i+1}-z_i\gg 0$ with $z_i-z_{i+1}\gg 0$ (note that $Z_{i,i+1}$ and $N_{i,i+1}$ are regular at $z_i=z_{i+1}$).

Therefore the analytic continuation of $F_0$ along the path $b_i$ amounts to the continuation of $ F_+^\circ(z_i,z_{i+1})$ along $b_i$, i.e., in the proper domain where both the above approximation of $F_0$ and $F_i$ are valid, we have in the asymptotic zone
\begin{eqnarray}\label{asyconstant}
F_0^{-1}\cdot F_i\sim N_{i,i+1}^{-1} F_+^\circ(z_i,z_{i+1})^{-1} F_-^\circ(z_i,z_{i+1}) N_{i,i+1}=N_{i,i+1}^{-1} S^{i,i+1} N_{i,i+1}=S^{i,i+1}. 
\end{eqnarray}
Here the first identity uses \eqref{reduce1} and \eqref{reduce2}, and the defining identity \eqref{defiSmatrix} of the Stokes matrix $S$ (the continuation of $ F_+^\circ(z_i,z_{i+1})$ along $b_i$ is same as the continuation of $Y_+(z)$ from $\IH_+$ to $\IH_-$ in a counterclockwise direction). The second identity follows from the fact that $N_{i,i+1}$ commutes with $S^{i,i+1}$. Since $F_0^{-1}\cdot F_i$ and $S^{i,i+1}$ are both constant, we get $F_0^{-1}\cdot F_i=S^{i,i+1}$. It proves that the monodromy of ${\rm gcKZ}$ equation along $b_i$ is given by $s_i(S(u))^{i,i+1}$.


Similarly, if we consider the asymptotic zone $z_{j+1}-z_j\gg 0, \ \forall j$ and $|\frac{z_1}{z_1-z_{i}}|\ll 1$ for any $i>1$, then the canonical solutions $F_0$ and $F_{-1}$ can be compared to the canonical solutions of equation \eqref{eq0}, i.e., $\kappa\frac{dF}{dz}=(u^{(1)}+\frac{\Omega_{01}+C_s^{(1)}}{z})\cdot F$. By the same argument, it follows that the monodromy along $\sigma$ is given by $\tau^{(1)}(K(u))$.
\qed

\vspace{2mm}
As a corollary, we have
\begin{thm}
For any $u$, the Stokes matrices $K(u)$ and $S(u)$ satisfy the Yang-Baxter and the $\tau$-twisted reflection equations.
\end{thm}
\pf It follows from Theorem \ref{monodromrep} and the braid relation in type $B$, as well as the fact $(\tau\otimes \tau)(S)=S.$
\qed

\begin{rmk}
Following Felder-Markov-Tarasov-Varchenko \cite[Theorem 3.1]{FMTV}, the gKZ equation has solutions taking the form of (confluent) hypergeometric integrals over twisted cycles (see e.g., \cite{Varchenko}). One can study the asymptotics of these solutions at the irregular singularity, and try to get an integral expression of the Stokes matrices. For example,
in the ${\rm gl}_2$ case, the gKZ equation, with two variables and valued in a two dimensional vector space, reduces to a confluent hypergeometric equation, and the integral solutions in \cite{FMTV} are related to the integral representations of confluent hypergeometric functions. The Stokes phenomenon of these solutions amounts to the different asymptotics of confluent hypergeometric function $_1F_1(a,b;z)$ as $z\rightarrow \infty$ from two different sectors. The comparison of the different asymptotics gives the Stokes matrices. More details can be found in the explicit computation of Stokes matrices in rank two \cite[Proposition 8]{BJL}. One should compare it to the ${\rm gl}_2$ example of KZ equation with three variables, and the hypergeometric function $_2F_1(a,b,c;z)$, see e.g., \cite[Section 1.1]{Varchenko}.
We refer the
reader to the book \cite{Varchenko} for a detailed exposition of the role of (multidimensional) hypergeometric functions in quantum algebras. Given this viewpoint, Theorem \ref{mainthm4} and Proposition \ref{twisteq} says that the "confluence" of quasi-triangular quasi-bialgebra structure becomes the quasi-triangular bialgebra structure (or equivalently, the "confluence" of associators become the twists). It is interesting to see the role of confluent hypergeometric functions in representation theory in more details and cases. Some more examples include studying the Stokes phenomenon of the affine KZ equation associated with Hecke type algebras (see \cite{Che3}), Dunkl operators, and the quantum KZ equation \cite{FR} (where in simplest case the q-confluent hypergeometric function will appear) with extra singularities.
\end{rmk}

\subsection{Variation of $u$}\label{section3}
In this subsection, we will study the derivations of the Stokes matrices $S(u)$ and $K(u)$ with respect to the parameter $u$, i.e., give a proof of Theorem \ref{mainthm2}. It needs the following lemma.
\begin{pro}\label{isodeformation}
The canonical solutions $F_k$, for any $k=-1,0,1...,n-1$, satisfy
\begin{eqnarray*}
\kappa d_\h F=\Big(\sum_{i=1}^n z_i du^{(i)}
+\sum_{\alpha\in \Delta}\frac{d\alpha}{\alpha}\Delta^{(n+1)}(C_{\frak k, \alpha})\Big)F-F\Big(\sum_{\alpha\in \Delta}\frac{d\alpha}{\alpha}(C_{\frak k,\alpha}^{(0)}+\sum_{i=1}^n C_{\alpha}^{(i)})\Big).
\end{eqnarray*}
Here $d_\h$ is the de Rham differential on $\h(\mathbb{R})$, $\Delta$ is the coproduct on $U(\g)$, and $\Delta^{(n)}:U(\g)\to U(\g)^{\otimes n}$
is the iterated coproduct recursively defined by
$\Delta^{(1)}=\id$, and $\Delta^{(n)}=(\Delta\otimes\id^{\otimes (n-2)})\circ
\Delta^{(n-1)}$ for $n\geq 2$.
\end{pro}
\pf Let us prove it for any fix $i\in\{-1,...,n-1\}$. Let us use the setting in Proposition \ref{solution}, i.e., consider the pull-back of the ${\rm gcKZ}$ equation under the map $P:\mathbb{C}^\times \times X_n \rightarrow X_n.$
Let us introduce the differential operators
\begin{eqnarray*}
&&D_z:=\kappa d_z-\Big(l\Big(\sum_{i=1}^n \xi_iu^{(i)}+ 
\frac{ 2\sum_{i<j}\Omega_{\frak k}^{ij}+\sum_{i}C_\frak k^{(i)}}{z}\Big)-r\Big(\sum_{i=1}^n \xi_iu^{(i)}+ 
\frac{ \sum_{i}C_0^{(i)}}{z}\Big)\Big)dz,\\
&&D_{u}:=\kappa d_u-\Big(l\Big(z\sum_{i=1}^n\xi_i du^{(i)}+\sum_{\alpha\in \Delta}\frac{d\alpha}{\alpha}\Delta^{(n+1)}C_{\frak k, \alpha}\Big)-r\Big(z\sum_{i=1}^n\xi_i du^{(i)}+\sum_{\alpha\in \Delta}\frac{d\alpha}{\alpha}(C_{\frak k,\alpha}^{(0)}+\sum_{i=1}^n C_{\alpha}^{(i)})\Big)\Big),
\end{eqnarray*}
where $l$ and $r$ denote the left and right multiplication.
First one checks that 
\[[D_{u},D_z]=0, \ and \ D_u^2=0, \ for \ any \ i, j=1,..., n.\]

Since the function $F_k$ is a solution, we get $D_zH_k=0$, see Proposition \ref{solution} for the definition of the function $H_k$. It implies that $D_zD_{u}H_k=D_{u}D_zH_k=0$. Thus to prove $D_{u}H_k=0$, by the uniqueness argument we only need to show that $D_{u}H_k$, as a function of $z$, tends to 0 as $z\rightarrow \infty$ in the sector $\wh{\IH}_+$. It can be seen as follows: plug the expansion $H_k\sim 1+ h_1z^{-1}+O(z^{-2})$ in $D_{u}H_k$, we see that the limit of $D_{u}H_k$ as $z\rightarrow\infty$ within $\wh{\IH}_+$ is given by
\begin{eqnarray*}
-[\sum_{i=1}^n\xi_idu^{(i)},h_1]-\sum_{\alpha\in\Delta}\frac{d\alpha}{\alpha}\Delta^{(n+1)}(C_{\frak k, \alpha})+\sum_{\alpha\in\Delta}\frac{d\alpha}{\alpha}(C_{\frak k, \alpha}^{(0)}+\sum_{i=1}^n C_{\alpha}^{(i)}).
\end{eqnarray*}
Then given the identity \eqref{simH} satisfied by $h_1$, one checks that the above expression vanishes. \qed

\vspace{2mm}
{\bf Proof of Theorem \ref{mainthm2}:}
From the definition of Stokes matrices and Proposition \ref{isodeformation}, 
Recall that Theorem \ref{monodromrep} (in particular the asymptotics in \eqref{asyconstant}) gives
\begin{eqnarray}\label{Sii+1}
F_0^{-1}\cdot F_i=S^{i,i+1}.
\end{eqnarray}
Following Proposition \ref{isodeformation}, we have
that both $F_0$ and $F_i$ satisfy $d_\h F=A F-B F$, where\[A=\sum_{i=1}^n z_i du^{(i)}
+\sum_{\alpha\in \Delta}\frac{d\alpha}{\alpha}\Delta^{(n+1)}(C_{\frak k, \alpha}) \hspace{2mm} \text{   and   } \hspace{2mm} B=\sum_{\alpha\in \Delta}\frac{d\alpha}{\alpha}(C_{\frak k,\alpha}^{(0)}+\sum_{i=1}^n C_{\alpha}^{(i)}).\]
Then taking the derivatives of both sides of \eqref{Sii+1} with respect to $u$ gives
\begin{eqnarray*}
d_\h S^{i,i+1}&=&-F_0^{-1}\cdot d_\h F_0\cdot  F_0^{-1} \cdot F_i+F_0^{-1}\cdot d_\h F_i\\
&=& -F_0^{-1}\cdot AF_i+BF_0^{-1}\cdot F_i+F_0^{-1} AF_i-F_0^{-1}\cdot F_iB\\
&=&[B, F_0^{-1}\cdot F_i]=[B, S^{i,i+1}].
\end{eqnarray*}
Taking respectively $i=1$ and $i=0$ in the above identity give rise to the expressions \eqref{isoeq1} and \eqref{isoeq1} in Theorem \ref{mainthm2}.\qed 

\vspace{2mm}
The equations \eqref{isoeq1} and \eqref{isoeq}, controlling the variation of $u$, are equivalent to the isomonodromy deformation equations of the differential equations \eqref{eq} and \eqref{eq0} respectively. We refer the reader to the theory of the isomonodromy deformation of linear systems of meromorphic ordinary differential equations by Jimbo, Miwa and Ueno \cite{JMU}. Theorem \ref{mainthm2} then states that the isomonodromy deformation amounts to gauge transformations of universal solutions of Yang-Baxter and reflection equations. In the end, we remark that it is interesting to study the asymptotics of the Stokes matrices $S(u,\kappa)$ and $K(u,\kappa)$ of the equations \eqref{eq} and \eqref{eq0}, as $\kappa\rightarrow 0$. Following \cite{Xu2}, the leading order in the asymptotics of $S(u,\kappa)$, as $\kappa\rightarrow 0$, is related to the theory of crystals. It is interesting to see if the leading order in the asymptotics of $K(u,\kappa)$ is related to the theory of canonical bases for quantum symmetric pairs introduced by Bao and Wang in \cite{BW}.

\begin{appendices}
\section{Canonical solutions, Stokes matrices and connection matrices}

Let us consider the linear system of meromorphic ordinary differential equations
\begin{eqnarray*}\label{nabla}
\frac{dF}{dz}=\Big(u+\frac{A}{z}\Big)\cdot F,
\end{eqnarray*}
where $F(z)\in {\rm GL}_n$, $u={\rm diag}(u_1,...,u_m)$ is a diagonal matrix, and $A\in {\rm gl}_n$. We divide $\{1,...,n\}$ into subsets $\{I_l\}_{l=1,...,k}$ such that $u_i=u_j$ if $i,j\in I_l$ for some $l$, and $u_i\ne u_j$ otherwise. We then assume that for any $l=1,...,k$, no two eigenvalues of the submatrix formed by the rows and columns from the index set $I_l$ differ by a non-zero integer.

The equation has an irregular singularity at $z=\infty$, and has a unique formal fundamental solution around $\infty$ taking the form (see e.g., \cite[Chapter 3]{Balser})
\begin{eqnarray*}\label{fomralsol}
\widehat{F}(z)=\widehat{H}(z)z^{[A]} e^{z u}, \ \ \ {\it for} \ \ \widehat{H}(z)=1+H_1z^{-1}+H_2z^{-2}+\cdot\cdot\cdot.
\end{eqnarray*}
Here $[A]$ takes the projection of $A$ to the centralizer of $u$ in ${\rm gl}_n$. In particular, if $u$ has distinct diagonal elements, $[A]$ takes the diagonal part of $A$. Although the radius of convergence of $\hat{H}(z)$ is in general zero, its Borel-Laplace transform (see \cite[Chapter 5]{Balser}) produces different holomorphic functions $H_i$, with the prescribed asymptotics $\hat{H}$ in certain different sectors $\wh{\rm Sect}_i$ of the complex plane, whose union covers a full neighborhood of the singularity $z=\infty$. See e.g., \cite[Chapter 8]{Balser}. In this way, one gets canonical fundamental holomorphic solutions $F_i(z)$ in each $\wh{\rm Sect}_i$. The mismatch of two sectoral solutions on the overlap of the corresponding sectors, known as the Stokes phenomenon, may be measured in terms of the transition matrix relating the two fundamental solutions. In what follows, we give more details and introduce some necessary notations for the present paper.

\begin{defi}\label{Stokesrays}
The {\it anti-Stokes rays} of the equation \eqref{nabla} are the directions
along which $e^{(u_i-u_j)z}$ decays most rapidly as $z\mapsto \infty$ for some $u_i\ne u_j$. The {\it Stokes sectors}
are the open regions of $\mathbb{C}$ bounded by two adjacent anti-Stokes rays.
\end{defi}

Given an initial Stokes sector ${\rm Sect}_0$, we label the anti-Stokes rays $d_1,d_2,...,d_{2l}$
going in a positive sense and starting on the positive edge of ${\rm Sect}_0$.
We denote by ${\rm Sect}_i$ the Stokes sector bounded by two adjacent anti-Stokes rays $d_i$ and $d_{i+1}$. Here indices are taken modulo $2l$, i.e., ${\rm Sect}_0={\rm Sect}(d_{2l}
, d_1)$. The following result can be found in e.g., \cite{Balser}\cite{BJL}\cite{MR}.
\begin{thm}\label{jurk}
On each ${\rm Sect}_i$, there is a unique (therefore canonical) holomorphic function $H_i:{\rm Sect}_i\to {\rm GL}_n$ such that the function
\[F_i(z)=H_i(z) e^{zu} z^{[A]}\]
satisfies equation \eqref{nabla}, and at the same time $H_i(z)$ can be analytically continued to a bigger sector $\widehat{{\rm Sect}_i}:=(d_i-\frac{\pi}{2},d_{i+1}+\frac{\pi}{2})$ and is asymptotic to $\widehat{H}$ as $z\rightarrow \infty$ within $\widehat{{\rm Sect}_i}$. 
\end{thm}

\begin{defi}
The {\it Stokes matrices} of the equation \eqref{nabla} (with respect to
to ${\rm Sect}_0$) are the matrices $S(A,u)$, $S_-(A,u)$ determined by
\[F_l=F_{0}\cdot e^{-\pi i [A]} S, \ \ \ \ \ 
F_{0}=F_l\cdot S_-e^{\pi i [A]}
\]
where the first (resp. second) identity is understood to hold in ${\rm Sect}_l$
(resp. ${\rm Sect}_0$) after $ F_0$ (resp. $ F_{l}$)
has been analytically continued counterclockwise. 
\end{defi}

Now if further no two eigenvalues of $A$ are differed by a non-zero integer, then the following fact is well-known (see e.g \cite[Chapter 2]{Wasow}).

\begin{lem}\label{le:nr dkz}
There is a unique holomorphic fundamental solution
$F_0(z)\in {\rm GL}_n$ of the system \eqref{nabla} on a neighbourhood of $\infty$ slit along the anti-Stokes ray $d_0$, such that $F_0\cdot z^{A}\rightarrow 1$ as $z\rightarrow 0$.
\end{lem}

\begin{defi}
The {\it connection 
matrix} $C(A,u)\in {\rm GL}_n$ of the system \eqref{nabla} (with respect to ${\rm Sect}_+$) is determined by 
$F_0(z)=F_+(z)\cdot C(u,A).$
\end{defi}

In a global picture, the connection matrix is related to the Stokes matrices
by the following monodromy relation, which follows from the fact that a simple negative loop (i.e., in cloclwise direction) around $0$ is a simple positive loop around $\infty$: 
\begin{eqnarray*}\label{monodromyrelation}
C(A,u)e^{2\pi i A}C(A,u)^{-1}=S_-(A,u)S(A,u).
\end{eqnarray*}

We remark that the Stokes matrices $S(A,u), S_-(A,u)$ will in general depend on the irregular data $u$ in \eqref{nabla}. Such dependence was studied by many authors, see e.g., \cite{JMU}.
\end{appendices}

\Addresses


\begin{thebibliography}{ZZZ}
\addtolength{\itemsep}{-1.5 em}
\setlength{\itemsep}{-5pt}
\bibitem{Ara}
S. Araki, {\it On root systems and an infinitesimal classification of irreducible symmetric spaces}, J.
Math. Osaka City Univ., 13(1), 1962.

\bibitem{BK} 
M. Balagovi${\rm \check{c}}$ and S. Kolb, {\it Universal K-matrix for quantum symmetric pairs}, J. reine angew. Math. 747 (2019), 299–353.

\bibitem{BW0}
H. Bao and W. Wang, {\it A new approach to Kazhdan-Lusztig theory of type B via quantum symmetric pairs}, Asterisque 402 (2018), vii+134pp.

\bibitem{BW} H. Bao and W. Wang, {\it Canonical bases arising from quantum symmetric pairs}, Invent. Math. 213 (2018), 1099–1177.

\bibitem{Balser}
W. Balser, {\it Formal power series and linear systems of meromorphic ordinary differential equations},
Springer-Verlag, New York, 2000.

\bibitem{BJL}
W. Balser, W.B. Jurkat and D.A. Lutz, {\it Birkhoff invariants and Stokes' multipliers for meromorphic linear differential equations}, J. Math. Anal. Appl. 71 (1979), 48-94.
 

\bibitem{Bri}
E. Brieskorn, {\it Die Fundamentalgruppe des Raumes der regularen Orbits einer endlichen
komplexen Spiegelungsgruppe}, Invent. Math. 12 (1971), 57–61.

\bibitem{Bro0}
A. Brochier, {\it A Kohno-Drinfeld theorem for the monodromy of cyclotomic KZ connections}, Comm. Math. Phys.
311 (2012), no. 1, 55–96.

\bibitem{Bro}
A. Brochier, {\it Cyclotomic associators and finite type invariants for tangles in the solid torus}, Algebr. Geom. Topol. 13 (2013), 3365–3409.


\bibitem{Ch}
L. Chekhov, {\it Teichm$\ddot{u}$ller theory of bordered surfaces}, SIGMA Symmetry and Integrability Geom. Methods Appl. 3 (2007), 066, 37 pp.
 
\bibitem{Che}
I. Cherednik, {\it Factorizing particles on a half-line and root systems}, Theoret. Math.
Phys 61 (1984), 977–983.


\bibitem{Che1}
I. Cherednik, {\it Generalized braid groups and local r-matrix systems}, Dokl. Akad. Nauk SSSR
307 (1989), no. 1, 49-53 (Russian); English translation in Soviet Math. Dokl. 40 (1990), no. 1, 43-48.

\bibitem{Che2}
I. Cherednik, {\it Calculation of the monodromy of some W-invariant local systems of type B, C and D}, Funktsional. Anal. i Prilozhen. 24 (1990), no. 1, 88–89. Translation in Funct. Anal. Appl. 24 (1990), no. 1, 78–79.

\bibitem{Che3}
I. Cherednik, {\it Lectures on Knizhnik-Zamolodchikov equations and
Hecke algebras}, Quantum many-body problems and representation theory, 1–96, MSJ
Memoirs 1, Math. Soc. Japan, Tokyo, 1998.



\bibitem{CNTY}
K. De Commer, S. Neshveyev, L. Tuset and M. Yamashita, {\it Ribbon Braided Module Categories, Quantum Symmetric Pairs and Knizhnik–Zamolodchikov Equations}. Commun. Math. Phys. 367, 717–769 (2019).

\bibitem{Drinfeld}
V. Drinfeld, {\it Quasi-Hopf algebras}, Algebra i Analiz 1 A989), no. 6, 114-148 (Russian);
English translation in Leningrad Math. J. 1 (1990), 1419-1457.

\bibitem{Drinfeld2}
V. Drinfeld, {\it On almost cocommutative Hopf algebras}, Leningrad Math. J. 1 (1990), no. 2,
321–342.

\bibitem{Dubrovin}
B. Dubrovin, {\it Geometry of 2D topological field theories}, Lecture Notes in Math, 1620 (1995).


\bibitem{En}
B. Enriquez, {\it Quasi-reflection algebras and cyclotomic associators,} Selecta Math. (N.S.) 13 (2007), no. 3, 391-463.


\bibitem{EE}
B. Enriquez and P. Etingof, {\it Quantization of classical dynamical r-matrices with nonabelian base,} Comm.
Math. Phys. 254 (2005), no. 3, 603–650,

\bibitem{EFK}
P. Etingof, I. Frenkel and A. Kirillov, {\it Lectures on representation theory and Knizhnik–Zamolodchikov equations}, Vol. 58 of Mathematical Surveys and Monographs (American Mathematical Society, Providence, RI, 1998). 

\bibitem{FMTV}
G. Felder, Y. Markov, V. Tarasov and A. Varchenko, {\it Differential equations compatible with KZ equations}, Math. Phys. Anal. Geom. 3 (2000), 139-177.

\bibitem{FR}
I. Frenkel and N. Reshetikhin, {\it Quantum affine algebras and holonomic difference equations}, Comm. Math. Phys., 146 (1), 1–60.

\bibitem{GL}
V. A. Golubeva and V. Leksin, {\it On a generalization of the Drinfeld-Kohno theorem}, Proceedings of the
Second ISAAC Congress, Vol. 2 (Fukuoka, 1999), 2000, pp. 1371–1386, Kluwer Acad. Publ., Dordrecht,

\bibitem{Jimbo}
M. Jimbo, {\it A $q$-difference analogue of $U\g$ and Yang-Baxter equation}, Lett.
in Math. Phys., 10 (1985), 63-69.

\bibitem{JMU}
M. Jimbo, T. Miwa and K. Ueno, {\it Monodromy preserving deformations of linear differential equations with rational coefficients $\uppercase\expandafter{\romannumeral1}$}, Physica 2D (1981), 306-352.

\bibitem{KZ}
V.G. Knizhnik and A.B. Zamolodchikov, {\it Current algebra and Wess-Zumino model in two dimensions}, Nuclear Phys. B 247 (1984), 83-103.


\bibitem{Kohno}
T. Kohno, {\it Monodromy representations of braid groups and Yang-Baxter equations},
Ann. Inst. Fourier 37 (1987), no. 4, 139-160.

\bibitem{Kolb}
S. Kolb, {\it Braided module categories via quantum symmetric pairs}, Proc. London Math. Soc.
(in press), arXiv:1705.04238v2 (2017).

\bibitem{KSS}
P.P. Kulish, R. Sasaki and C. Schwiebert, {\it Constant solutions of reflection equations
and quantum groups}, J. Math. Phys. 34 (1993), no. 1, 286–304.

\bibitem{KSly}
P.P. Kulish and E. Sklyanin, Algebraic structures related to reflection equations, J.
Phys. A: Math. Gen. 25 (1992), 5963–5975.

\bibitem{Lei}
A. Leibman, {\it Some monodromy representations of generalized braid groups}, Comm. Math. Phys. 164 (1994),
no. 2, 293–304.

\bibitem{Let}
G. Letzter, {\it Symmetric pairs for quantized enveloping algebras}, J. Algebra 220 (1999),
729–767.

\bibitem{Let2}
G. Letzter, {\it Coideal subalgebras and quantum symmetric pairs}, New directions in Hopf
algebras (Cambridge), MSRI publications, vol. 43, Cambridge Univ. Press, 2002,
pp. 117–166.

\bibitem{LW}
J.-H. Lu and A. Weinstein, {\em Poisson Lie groups, dressing transformations, and Bruhat decompositions}, J. Diff. Geom. 31 (1990), 501 - 526.

\bibitem{MR}
B. Malgrange and J.-P. Ramis, {\it Fonctions multisommables}, Ann. Inst. Fourier (Grenoble) 42 (1992),
no. 1-2, 353–368.

\bibitem{MolRag}
A. I. Molev and E. Ragoucy,
{\it Symmetries and invariants of twisted quantum algebras and associated Poisson algebras}, Rev. Math. Phys. 20 (2008), 173-198.

\bibitem{Noumi}
M. Noumi, {\it Macdonald’s Symmetric Polynomials as Zonal Spherical Functions on Some Quantum Homogeneous Spaces}, Adv. Math. 123 (1996), 16–77.

\bibitem{NDS} M. Noumi, M.S. Dijkhuizen, and T. Sugitani, Multivariable Askey-Wilson polynomials
and quantum complex Grassmannians, AMS Fields Inst. Commun. 14 (1997), 167–177.

\bibitem{NS} M. Noumi and T. Sugitani, {\it Quantum symmetric spaces and related q-orthogonal polynomials,} Group theoretical methods in physics (Singapore) (A. Arima et. al., ed.),
World Scientific, 1995, pp. 28–40.

\bibitem{Res2} 
N. Reshetikhin, {\it The Knizhnik-Zamolodchikov system as a deformation of the isomonodromy problem},
Lett. Math. Phys. 26 (1992), no. 3, 167–177.

\bibitem{Skl}
E.K. Sklyanin, {\it Boundary conditions for integrable quantum systems}, J. Phys. A 21
(1988), 2375–2389.


\bibitem{TL} V. Toledano Laredo, {\it Quasi--Coxeter quasitriangular quasibialgebras and the Casimir connection}, arXiv:1601.04076.


\bibitem{TLXu}
V. Toledano Laredo and X. Xu, {\it Stokes phenomenon, Poisson Lie groups and quantum groups.} In preparation.

\bibitem{tD}
T. tom Dieck, {\it Categories of rooted cylinder ribbons and their representations}, J. reine
angew. Math. 494 (1998), 36–63.

\bibitem{tDH}
T. tom Dieck and R. Haring-Oldenburg, {\it Quantum groups and cylinder braiding}, Forum
Math. 10 (1998), no. 5, 619–639.

\bibitem{TK}
A. Tsuchiya and Y. Kanie, {\it Vertex operators in two dimensional conformal field theory on $\mathbb{P}^1$
and monodromy representation of braid groups}, Adv. Stud, in Pure Math., 16 (1988),
297-372.


\bibitem{Varchenko}
A. Varchenko, {\it Hypergeometric functions and representation theory of Lie algebras and quantum groups}, Advanced Series in Mathematical Physics, Vol. 21, World Scientific (1995).

\bibitem{Wasow}
W. Wasow, {\em Asymptotic expansions for ordinary differential equations}, Wiley Interscience, New York,
1976.

\bibitem{Xu}
X. Xu, {\it Stokes phenomenon and Yang-Baxter equations}, Commun. Math. Phys. 377 (2020), 149-159.


\bibitem{Xu2}
X. Xu, {\em Closure of Stokes matrices I: caterpillar points and Alekseev-Meinrenken diffeomorphisms}, arxiv:1912.07196.
\end{thebibliography}
\end{document}